\documentclass{amsart}
\usepackage{amsmath}
\usepackage{graphicx}
\usepackage{amsfonts}
\usepackage{amssymb}

\newtheorem{theorem}{Theorem}
\theoremstyle{plain}

\newtheorem{corollary}{Corollary}

\newtheorem{lemma}{Lemma}

\newtheorem{proposition}{Proposition}
\newtheorem{remark}{Remark}

\numberwithin{equation}{section}

\catcode`@=11
\def\plainvspace@{\def\vspace##1{\noalign{\vskip##1}}}
\def\plainLet@{\relax\iffalse{\fi\let\\=\cr\iffalse}\fi}
\def\aligned{\,\vcenter\bgroup\plainvspace@\plainLet@\openup\jot\m@th\ialign
  \bgroup \strut\hfil$\displaystyle{##}$&$\displaystyle{{}##}$\hfil\crcr}
\def\endaligned{\crcr\egroup\egroup}
\def\matrix{\,\vcenter\bgroup\plainLet@\plainvspace@
    \normalbaselines
  \m@th\ialign\bgroup\hfil$##$\hfil&&\quad\hfil$##$\hfil\crcr
    \mathstrut\crcr\noalign{\kern-\baselineskip}}
\def\endmatrix{\crcr\mathstrut\crcr\noalign{\kern-\baselineskip}\egroup
                \egroup\,}

\def\cases{\left\{\,\vcenter\bgroup\plainvspace@
     \normalbaselines\openup\jot\m@th
      \plainLet@\ialign\bgroup$\displaystyle{##}$\hfil&\quad$\displaystyle{{}##}$\hfil\crcr
      \mathstrut\crcr\noalign{\kern-\baselineskip}}

\begin{document}
\title{A sharp   inequality and its applications}
\author  {Suyu Li and Meijun Zhu}
\address{Department of Mathematics\\
The University of Oklahoma\\
Norman, OK 73019\\ }
\begin{abstract}
We establish an analog Hardy inequality with sharp constant
involving exponential weight function. The special case of this
inequality (for $n=2$)  leads to a direct proof of Onofri
inequality on $S^2$.

\end{abstract}








 \maketitle

\section{Introduction}
The classical Hardy inequality says that for any non-negative
function $f(x)$ on $ [0, +\infty)$, if $F(x)=\int_0^x f(t) dt$,
then
$$
\int_0^\infty (\frac Fx)^k dx \le (\frac k{k-1})^k\int_0^\infty
f^k dx,$$ where $k>1$ is a given parameter. See, for example,
Inequality 327 in  the book by  Hardy, Littlewood and Polya
\cite{HLP}. It is important to note that the constant $( \frac
k{k-1})^k$ is the optimal one. Using H\"older inequality Hardy and
Littlewood  were able to derive that
$$\int_0^\infty \frac{F^l}{x^{l-\alpha}}dx \le (\frac k{k-1})^k
(\int_0^\infty f^k dx)^{\frac l k},$$ where $l\ge k$ and
$\alpha=l/k-1$. It was quite clear to them that the constant is
not optimal for $l>k$. Though they guessed what is the best
constant, it was later  proved by Bliss, who obtained  nowadays
the famous Bliss Lemma (see the interesting papers \cite{HL} and
\cite{BL}):

\medskip

\noindent {\bf Bliss Lemma}: {\it Let $k,\ l$ be constants, such
that $l>k>1$, and let $f(x)$ be a   non-negative measurable
function in the intervals $0 \le x<\infty$, such that the integral
$J=\int_0^\infty f^k dx$ is finite. Then the integral $y=\int_0^x
fdx$ is finite for every $x$ and \begin{equation} I=\int_0^\infty
\frac{y^l}{x^{l-\alpha}} dx \le C_{b} J^{l/k},\label{bl-1}
\end{equation}
where $$\alpha=\frac lk -1, \ \ \ \ C_{b}=\frac
1{l-\alpha-1}[\frac{\alpha \Gamma(l/\alpha)}{\Gamma(1/\alpha)
\Gamma((l-1)/\alpha)} ]^\alpha.$$ The equality in (\ref{bl-1})
holds if and only if $f(x)=c/(1+d x^\alpha)^{(\alpha+1)/\alpha}$
for some positive constants $c,\ d$. }

\medskip

Bliss Lemma  later (after more than forty years) became a crucial
ingredient in the proof of  sharp Sobolev inequality by Aubin
\cite{AU}, and Talenti \cite{TA} respectively. The latter
inequality has played essential role in the resolution of the
Yamabe problem, which mainly concerns about finding a canonical
metrics with constant scalar curvature on compact manifolds with
dimension higher than or equal to three (see the geometric and
analytic forms of sharp Sobolev inequalities in the appendix).


 The Yamabe
problem can also be viewed as the higher dimensional analogue to
the uniformization theorem for two dimensional manifolds. The
analytic approach to the re-proof of the uniformization theorem
seems to be initiated  by Berger \cite{BE}.  While for the Yamabe
problem for manifold with positive Yamabe constant, one seeks a
new metric in the same conformal class (with fixed volume)  which
yields the smallest total scalar curvature, in the analytic
approach to the re-proof of the uniformization theorem on
topological spheres one looks for  a new metric in the same
conformal class (with fixed area) which has the smallest Liouville
energy, see, for example, Hamilton \cite{HA}, Chow \cite{CHO}, or
the recent paper by Chen and Zhu \cite{CZ}. The core inequality in
such an argument is the Onofri inequality (see the precise form in
the appendix). Recently we showed in \cite{LZ} that one can derive
the Onofri inequality directly from Trudinger's inequality.
Comparing the proof of sharp Sobolev inequality with that of
Onofri inequality, we feel that there is an undiscovered calculus
inequalities, which turns out to be the main theorem of this
paper.


\begin{theorem}{\label{main}}
1). Let $n>1$ be given.  For any nonnegative function $u\in C^1[0,
+\infty)$ with $u(0)=0$,
\begin{equation}
\ln \int_0^{+\infty} \frac {e^{nu}}{e^{nr}}dr \le
(\frac{n-1}n)^{n-1} \int_0^{+\infty}|u_r|^n dr +C_n, \label
{hardyet}
\end{equation}
  where the constant
 $C_n$ is given by
$$C_n= \int_{0}^1 \frac 1{t}(1-
 (1-t)^{n-[n]})
\ d t+ \sum_{i=1}^{[n]-1} \frac {1}{(n-i)},$$ $[n]$ is the integer
part of $n$. Both constants $(\frac{n-1}n)^{n-1}$ and $C_n$ are
optimal, and the equality never holds.

2). For any nonnegative function $u\in C^1[0, +\infty)$ with
$u(0)=0$,
\begin{equation}
\ln \int_0^{+\infty} \frac {e^{u}}{e^{r}}dr \le
 \int_0^{+\infty}|u_r| dr. \label {n=1}
\end{equation}
\end{theorem}

\medskip

We first prove the above  inequality with a larger coefficient in
Section 2 (Proposition \ref{rough-1} below). The argument  is
elementary and simple. It needs to be pointed out that for $n>1$
being an integer, Theorem \ref{main} can be read out from Theorem
1.3 in \cite{LZ}. For general positive constant, it seems
impossible to prove Theorem \ref{main} from that theorem, rather,
Theorem \ref{main} provides an alternative proof of that theorem
(Corollary \ref{LZcor} in this paper). Recall the original proof
of Theorem 1.3 in \cite{LZ} does rely on Trudinger's inequality.
Quite interestingly, we also recall that Moser \cite{MO} used a
similar argument to give a very simple proof of the improved
Trundinger's inequality (with best constant):
\begin{corollary} ({\bf Weak Moser's inequality})
Let $\Omega \subset {\bf R}^n$ (for $n \ge 2$) be a smooth bounded
domain. For any $\beta <n \omega_{n-1}^{1/(n-1)}$, there is a
constant $C(\Omega,\beta)$ depending on the volume of $\Omega$ and
$\beta$, such that for all $u\in W^{1, n}_0(\Omega)$ with
$\int_\Omega |\nabla u|^n dx \le 1$,
$$\int_\Omega e^{\beta u^{\frac n{n-1}}} dx \le C(\Omega,\beta).$$
\label{TM-1}
\end{corollary}
Here and throughout this paper, we use $\omega_n$ for the volume
of unit sphere $S^n$ in $R^{n+1}$.  This result is slightly weaker
than Moser's inequality since it does not include the case of
$\beta= n \omega_{n-1}^{1/(n-1)}$. It seems that one needs the
argument due to Moser \cite{MO}, or Carleson and Chang \cite{CC}
to cover this extremal case.

In Section 3 we will show how to improve the rough inequality
(Proposition \ref{rough-1}) and complete the proof of  the main
theorem.  One particular reason that we can achieve this (but not
for Moser's inequality) is that we can classify all extremal
functions.

 As the
Bliss Lemma yields sharp Sobolev inequality, in Section 4 we will
show that  Theorem \ref{main} can be used to give a more direct
proof of the
 Onofri inequality (thus without even using Trudinger's inequality).
In fact, let $B_r(0) \subset R^n$ (now $n$ is an integer greater
than or equal to two) be a
 ball in $R^n$ with radius $r$ centered at the origin, and
$$
D^b_{a}(B_r(0))=\{f(y)\ : \ f(y)-b\in W_0^{1,n}(B_r(0)), \   \
\int_{B_r(0)} e^{nf}dy=a\},
$$
where $a $ is  a constant  satisfying $ a> \frac {\omega_{n-1} r^n
e^{nb}}n$. We will show that Theorem \ref{main} yields

\begin{corollary} {\bf (Local sharp inequality for $n=2$)}
$$\inf_{w \in D^b_{a}(B_r)}  \int_{B_r}|\nabla w|^2 d y= 4
\pi \cdot (\ln \frac{a e^{-2b} }{\pi r^2} +\frac {\pi r^2}{a
e^{-2b} }-1).
$$
 \label{mainthm-1}
\label{LZcor1}
\end{corollary}
It is known now that this corollary implies
 Onofri inequality on $S^2$, see, Li and Zhu
\cite{LZ}. For readers' convenience we include a complete proof of
the Onofri inequality in Section 4.

 In Section 4
we shall also discuss the applications of the main theorem to
other  geometric problems.  For readers' convenience, we present
both geometric and analytic forms of sharp Sobolev inequality on
$S^n$ (for $n \ge 3$) and Onofri inequality on $S^2$ in the
appendix.

\medskip
\noindent{\bf ACKNOWLEDGMENT. } The work of M. Zhu is partially
supported by the NSF grant DMS-0604169.

\section{Rough inequality}
We shall establish two elementary calculus inequalities in this
section. The first one will be used to prove the main theorem, and
the second one will be used to derive Corollary \ref{TM-1}.

\begin{proposition}
(1). Let $n> 1$  and $\beta_{0}>(\frac{n-1}n)^{\frac{n-1}n}$.
There is a constant $c_1(\beta_{0})$, such that for any $u(r) \in
C^1[0, +\infty)$ satisfying $ u(0)=0$,
\begin{equation}
\ln \int_0^\infty e^{n(u-r)} dr \le \beta_{0}^n
\int_0^\infty|u_r|^n dr+c_1(\beta_{0}).
\end{equation}

(2). For $u(r) \in C^1[0, +\infty)$ satisfying $ u(0)=0$,
$$
\ln \int_0^\infty \frac{e^u}{e^r} dr \le \int_0^\infty|u_r| dr.$$
\label{rough-1}
\end{proposition}

\medskip

 For $n=1$, the above  is an optimal inequality. For $n \ge 2$
we will improve the inequality by variational method in next
section.

\begin{proof}
Let  $u(r) $ be any  function in $C^1[0, +\infty)$ satisfying $
u(0)=0$.  We have
$$u(r) \le \int_0^\infty |u_r| dr,$$
thus
$$
\int_0^\infty \frac{e^u}{e^r} dr \le \exp\{\int_0^\infty|u_r|
dr\},$$ which yields
$$
\ln \int_0^\infty \frac{e^u}{e^r} dr \le \int_0^\infty|u_r| dr.$$

Now, for  given $n>1$ and positive parameter $\beta >0$, we have
\begin{align*}
u(r)& =\int_0^r u_r dr \le (\int_0^r|u_r|^n dr)^{1/n} \cdot
r^{\frac{n-1}n}\\
&\le \frac {\beta^n \int_0^r|u_r|^n dr}{n}+\frac{\beta^{-\frac
n{n-1}} r}{n/(n-1)}.
\end{align*}
Thus
\begin{align}\label{rough}
\int_0^\infty \frac{e^{nu}}{e^{nr}} dr &\le \int_0^\infty
\frac{\exp\{ \beta^n \int_0^\infty|u_r|^n dr+(n-1)\beta^{-\frac
n{n-1}}r\}}{e^{nr}} dr\\
&=\exp \{ \beta^n \int_0^\infty|u_r|^n dr\} \cdot \int_0^\infty
e^{[(n-1)\beta^{-\frac n{n-1}}-n]r} dr. \nonumber
\end{align}
If we choose
\begin{equation}\label{2-1}
\beta=\beta_{0}>(\frac{n-1}n)^{\frac{n-1}n},
\end{equation}
then
$$
\int_0^\infty e^{[(n-1)\beta_{0}^{-\frac n{n-1}}-n]r}
dr=c(\beta_{0})$$ is a finite number depending  on $\beta_{0}$. It
follows that
$$\ln \int_0^\infty e^{n(u-r)} dr \le  \beta_{0}^n \int_0^\infty|u_r|^n
dr+c_1(\beta_{0})
$$
for  $c_1(\beta_{0})=\ln c(\beta_{0})$.

\end{proof}

It is obvious in the above proof that $c(\beta_{0}), c_1(\beta_{0})
\to +\infty$ as $\beta_{0}\to (\frac{n-1}n)^{\frac{n-1}n}.$ We need
another argument to derive the main theorem.

\begin{remark}\label{rem2-1}
From (\ref{rough}) we can see that for $\beta_{0}$ satisfying
(\ref{2-1}),
\begin{align}\label{rem2}
\int_R^\infty \frac{e^{nu}}{e^{nr}} dr &\le \exp \{ \beta_{0}^n
\int_0^\infty|u_r|^n dr\} \cdot \int_R^\infty
e^{[(n-1)\beta_{0}^{-\frac n{n-1}}-n]r} dr\\
&=o_R(1)\exp \{ \beta_{0}^n \int_0^\infty|u_r|^n dr\}, \nonumber
\end{align}
where $o_R(1)\to 0$ as $R\to \infty$.
\end{remark}

We now compare this with Moser's proof of Trudinger's inequality
\begin{lemma} For $n>1$, $a>0$  and $\beta<n a^{\frac 1{1-n}}$, there is a constant
$C_{\beta,a}$ depending only on $\beta$ and $a$, such that for any
nonnegative function $u \in C^1[0, +\infty)$ with $u(0)=0$ and
$\int_0^\infty|u_r|^n dr\le a$,
$$\int_0^\infty \frac{e^{\beta u^{\frac n{n-1}}}}{e^{nr}} dr \le
C_{\beta,a}.
$$ \label{propM-1}
\end{lemma}
\begin{proof}
For given $n>1$  we have
$$
u(r) =\int_0^r u_r dr \le (\int_0^r|u_r|^n dr)^{1/n} \cdot
r^{\frac{n-1}n} \le  a^{\frac 1n} r^{\frac{n-1}n} . $$ Thus for
any positive parameter $\tau
>0$,
\begin{equation}\label{roughM-1}
\int_0^\infty \frac{e^{\tau  u^{\frac n{n-1}}}}{e^{nr}} dr \le
\int_0^\infty \exp\{\tau a^{\frac 1{n-1}}-n\}r  dr. \end{equation}
The right hand side of the above inequality is bounded if we
choose $\tau=\beta<n a^{\frac 1{1-n}}$.
\end{proof}

Based on Lemma \ref{propM-1}, one can verify Corollary \ref{TM-1}
as follows.

Due to the rearrangement and rescaling, we only need to prove
Corollary \ref{TM-1} when $\Omega=B_1(0)$ and $u\in C_0^1(B_1(0))$
is radially symmetric and nonnegative.

From $\int_{B_1} |\nabla u|^n dx\le 1$, we know that (let $r=-\ln
s$)
\begin{align*}
1 \ge \int_{B_1}  |\nabla u|^n dx =\omega_{n-1} \int_0^1 |u_s|^n
s^{n-1} ds= \omega_{n-1} \int_0^\infty |u_r|^n dr.
\end{align*}
Also,
\begin{align*}
\int_{B_1} e^{\beta u^{\frac n{n-1}}} dx=\omega_{n-1} \int_0^1
e^{\beta u^{\frac n{n-1}}} s^{n-1} ds=\omega_{n-1} \int_0^\infty
\frac{e^{\beta u^{\frac n{n-1}}}}{e^{nr}} dr.
\end{align*}
One immediately has Corollary \ref{TM-1} by  using Lemma
\ref{propM-1} with $a=\omega_{n-1}^{-1}.$

\section{Sharp Inequality}
We shall prove the main theorem in this section. Since the case of
$n=1$ has been settled by Proposition \ref{rough-1}, we will focus
on the case of $n>1$. For given $a>0$, define
\begin{equation}\label{set-1}
D_a^n:=\{ u(r)\in W^{1,n}(R^+) \ : \ u(0)=0, \ \int_0^\infty
\exp\{nu-nr\} dr=a\}.
\end{equation}

\begin{lemma}\label{sharp}
There is a $v \in D_a^n$ such that
$$
\int_0^\infty |v_r|^n dr =\inf_{u \in D_a^n} \int_0^\infty |u_r|^n
dr:=\xi. $$
\end{lemma}

\begin{proof}
Let $\{v^i\}$ be a  minimizing sequence of $\inf_{u \in D_a^n}
\int_0^\infty |u_r|^n dr$. Then
$$v^i \rightharpoonup v \ \ \ \mbox{in} \ W^{1,n}(R^+), \ \ \mbox{and}
\ \ \ \int_0^\infty |v_r|^n dr \le \underline{\lim}_{i\to \infty}
 \int_0^\infty |v^i_r|^n dr=\xi $$ for some $v\in W^{1,n}(R^+).$
We  need to verify $v \in D_a^n.$

First, from (\ref{rem2}), we know that for $w= v^i,$ or $v$:
$$\int_R^\infty \frac{e^{nw}}{e^{nr}} dr =o_R(1).$$ On the other
hand,  it follows from  the embedding $H^1(0, R)\hookrightarrow
C^{0, 1/2}(0, R)$ and Arzela-Ascoli lemma  that
$$
\lim_{ i \to \infty} \int_0^R \exp\{nv^i -nr\} dr = \int_0^R
\exp\{nv -nr\} dr.$$ Letting $i, \ R\to \infty$, we have
$\int_0^\infty \exp\{nv -nr\} dr=a$, that is $v \in D_a^n.$
\end{proof}

We now begin the proof of the main theorem.
\begin{proof}We only need to consider nontrivial nonnegative
functions. For $a>1/n$, let $v$ be the minimizer of $\inf_{u \in
D_a^n}\int_0^\infty |u_r|^n dr$. It is easy to see that $v_r \ge
0$.  So it satisfies the following Euler-Lagrange equation:
\begin{equation}\label{equ-1}
v_r^{n-2} v_{rr}=-\tau e^{nv -nr}, \ \ \ \ v(0)=0 \end{equation}
for some $\tau>0$. Though it is not obvious how to obtain the
general solution from the uniqueness of the ordinary differential
equation {since $v_r$ could be zero},  one can follow the argument
given by Carleson and Chang (\cite{CC}, page 123) to show that the
general solution to (\ref{equ-1}) is given by
$$ v(r)=\ln \frac {1}{\lambda_0 + e^{-nr/(n-1)}}-\frac 1n \ln \frac{\tau }
{ (\frac{n}{n-1})^n \lambda_0} ,$$
where $ \lambda_0 $ is a positive constant and
$\tau=\frac{(\frac{n}{n-1})^n\lambda_0}{(\lambda_0+1)^n}. $ Thus
\begin{equation}\label{sol}  v(r)=\ln \frac
{\lambda_0+1}{\lambda_0 + e^{-nr/(n-1)}}. \end{equation} Since
$a=\int_0^\infty e^{nv-nr} \ d r ,$ we have
\begin{align*}
a&=\int_0^\infty\left(\frac{\lambda_0+1}{\lambda_0+e^{-nr/(n-1)}}\right)^n
e^{-nr}\ d r\\
&=\int_0^1 \left(\frac{\lambda_0+1}{\lambda_0+s^{n/(n-1)}}\right)^n
s^n (\frac1{s}) \ d s \\
&=(\lambda_0+1)^n  \int_0^1
\frac{s^{n-1}}{(\lambda_0+s^{n/(n-1)})^n} \ d s\\
&=(\lambda_0+1)^n \frac {s^n}{n\lambda_0 (\lambda_0 +
s^{n/(n-1)})^{n-1}} \bigr|^{1}_{s=0}\\
&=\frac{\lambda_0+1}{n\lambda_0} ,\\
\end{align*}
That is
\begin{equation}\label{a-1} \lambda_0=\frac1{na-1} .
\end{equation}
We compute
\begin{align*}
\int_0^\infty |v_r|^n dr &= \int_0^\infty \left| \frac{\frac{n}{n-1}
e^{-nr/(n-1)}}{\lambda_0+e^{-nr/(n-1)}}\right|^n \ d r \\
&= (\frac{n}{n-1})^n \int_0^\infty
\left(\frac{e^{-nr/(n-1)}}{\lambda_0+e^{-nr/(n-1)}}\right)^n \ d r\\
&=(\frac{n}{n-1})^{n-1} \int_0^{1/\lambda_0} \frac{t^{n-1}}{(1+t)^n} \ d t\\
&=(\frac{n}{n-1})^{n-1} \int_1^{1/\lambda_0+1}
\frac{(\tau-1)^{n-1}}{\tau^n} \ d \tau\\
&=-(\frac{n}{n-1})^{n-1} \int_1^{1/\lambda_0+1}(\tau -1)(1-\frac
1 {\tau})^{n-2} \ d \frac 1{\tau}\\
& =(\frac{n}{n-1})^{n-1} \int_{{\lambda_0}/(\lambda_0+1)}^1 (\frac
1 t-1)(1-t)^{n-2} \ d t.
\end{align*}
Using (\ref{a-1}) we have: If $n \in  $ {\bf N}$,$
\begin{align*}
 \int_0^\infty |v_r|^n \ d r
&=(\frac{n}{n-1})^{n-1} \cdot
\{ -\int_{{\lambda_0}/(\lambda_0+1)}^1 (1-t)^{n-2} \ d t \\
&- \cdots -
 \int_{{\lambda_0}/(\lambda_0+1)}^1 (1-t) \ d t - \int_{{\lambda_0}/(\lambda_0+1)}^1 (1-\frac 1 {t}) \ d t
 \}\\
&=(\frac{n}{n-1})^{n-1} \{ \ln \frac {\lambda_0+1}{\lambda_0}
-\frac1{\lambda_0+1}-\sum_{i=1}^{n-2} \frac {(\frac
1{\lambda_0+1})^{n-i}}{n-i} \} \\
&=(\frac{n}{n-1})^{n-1} \{\ln (na) - \sum_{i=1}^{n-1} \frac
{(na-1)^{n-i}}{(n-i)(na)^{n-i}} \};
\end{align*}
For general $n>1$, we have
\begin{align*}
 \int_0^\infty |v_r|^n \ d r &=(\frac{n}{n-1})^{n-1} \cdot
\{ -\int_{{\lambda_0}/(\lambda_0+1)}^1 (1-t)^{n-2} \ d t \\
&- \cdots -
 \int_{{\lambda_0}/(\lambda_0+1)}^1 (1-t)^{n-[n]} \ d t + \int_{{\lambda_0}/(\lambda_0+1)}^1
 \frac 1{t}(1-t)^{n-[n]} \ d t
 \}\\
&=( \frac n{n-1} )^{n-1} \{\int_{\lambda_0/(\lambda_0+1)}^1 \frac
1{t}(1-t)^{n-[n]} \ d t - \sum_{i=1}^{[n]-1} \frac {( \frac
1{\lambda_0+1})^{n-i}}{n-i} \} \\
&= ( \frac n{n-1} )^{n-1} \{\int_{\lambda_0/(\lambda_0+1)}^1 \frac
1{t}+ \int_{\lambda_0/(\lambda_0+1)}^1 \frac 1{t}((1-t)^{n-[n]}-1)
\ d t - \sum_{i=1}^{[n]-1} \frac {( \frac
1{\lambda_0+1})^{n-i}}{n-i} \}\\
 &\geq (\frac{n}{n-1})^{n-1} \{\ln (na) -  \int_{\lambda_0/(\lambda_0+1)}^1 \frac 1{t}(1-
 (1-t)^{n-[n]})
\ d t- \sum_{i=1}^{[n]-1} \frac {(na-1)^{n-i}}{(n-i)(na)^{n-i}} \}
,
\end{align*}
where $[n]$ is the integer part of $n$. Let $a \to \infty$, then
$\lambda_0 \to 0$ by (\ref{a-1}).  We know that $C_n$ is optimal.
The proof  is completed.
\end{proof}

\begin{remark}For negative function $u$, we can certainly improve the
inequalities. In particular, similar argument will yield Theorem
1.3 (ii) in \cite{LZ} for integer $n>1$. Since we do not have
meaningful applications for this inequality so far, we shall skip
 details here.
\end{remark}


\section{Applications}
We shall show in this section that Theorem \ref{main} implies
Corollary \ref{LZcor1}, as well as Theorem 1.3 in \cite{LZ}, and
then shall derive th e Onofri inequality from Corollary
\ref{LZcor1}.

We first prove Corollary \ref{LZcor1}. Let $v\in D^2_\alpha$
(recalling the notation in (\ref{set-1}). We have, from the proof
of Theorem \ref{main}, that
\begin{equation}
\inf_{v \in D_\alpha^2} \int_0^\infty |v_r|^2 \ d r \geq 2 \{\ln
(2\alpha) +\frac{1}{2 \alpha } -1 \}, \label{n=2}
\end{equation}
where $\int_0^\infty e^{2v-2r}dr=\alpha$. For $w\in
D_a^0(B_1(0))$,
$$
\int_{B_1}  |\nabla w|^2 dx =2\pi \int_0^1 |w_s|^n s ds= 2 \pi
\int_0^\infty |w_r|^2 dr.
$$
and
$$
\int_{B_1} e^{2w } dx=2 \pi \int_0^1 e^{2w} s ds=2 \pi
\int_0^\infty \frac{e^{2 w}}{e^{2r}} dr.
$$
Combing with (\ref{n=2}), we have
$$
\inf_{w \in D^0_{a}(B_r)} \int_{B_1}  |\nabla w|^2 dx = 4 \pi
\cdot (\ln \frac{\int_{B_1} e^{2w } dx}{\pi} +\frac
{\pi}{\int_{B_1} e^{2w } dx}-1).$$ After rescaling and shifting,
we get Corollary \ref{LZcor1}.

\medskip

In the same spirit, we easily obtain
\begin{corollary}Let $u\in C^1(B_1)$ be a nonnegative
function satisfying $u=0$ on $\partial B_1$
$$
\ln \frac{n \int_{B_1} e^{nu}}{\omega_{n-1}} < (\frac{n-1}n)^{n-1}
\omega_n^{-1} \int_{B_1} |\nabla u|^n+F(1).$$
 where
$$F(1) = 1+\frac 12 +\cdots+\frac 1{n-1}.$$ \label{LZcor}
\end{corollary}

The fact that the strict inequality holds on a bounded domain
coincides with the one that the strict sharp Sobolev inequality
holds on a bounded domain. Corollary \ref{LZcor} was first proved
in \cite{LZ} using Trudinger's inequality. The proof presented
here does not rely on Trudinger's inequality. Inequality in
Corollary \ref{LZcor} was refereed as {\it local sharp inequality}
in \cite{LZ}, which is  easily adapted  for manifolds. See related
topics in Chen and Zhu \cite{CZ}.

\medskip

Finally, we shall show that one can prove the Onofri inequality
(see both forms of the inequality in appendix) using Corollary
\ref{LZcor1}.

Due to the rearrangement, we only need to prove Onofri inequality
for $u\in C^1(S^2)$ which
 depends only on $x_3$  and
is monotonically decreasing in $x_3$. Also, we can assume that
$u(x_3)\mid_{x_3=1}=0$ (otherwise, we replace $u(x)$ by
$u(x)-u(1))$. We can approximate $u$ by  a sequence of functions
$u_i \in C^1(S^2)$ such that $u_i(x)=u_i(x_3)$ is monotonically
decreasing in $x_3$, and $u_i(x)=0$ in the geodesic ball
$B_{1/i}(N)$ of the north pole $N$ for $i\in {\bf N}$. Denote
$S^2_i:=S^2 \setminus B_{1/i}(N)$.

Let $\Phi$: $x \in S^2 \to y \in R^2$ be a stereographic
projection given by
$$
x_i=\frac {2 y_i}{1+|y|^2}, \ \ \ \ \ \ for \ \ \ i=1, 2;
$$
and
$$
x_{3} =\frac{|y|^2-1}{|y|^2+1}.
$$
Denote
$$
g_0=\sum_{i=1}^{3}dx_i^2=(\frac{2}{1+|y|^2})^2 dy^2:=e^{2
\varphi(y)} d y^2.
$$
Thus
$$\varphi(y)=\ln\frac  2{1+|y|^2}.
$$
It is easy to check that $\varphi(y)$ satisfies
\begin{equation}
-\Delta \varphi=e^{2\varphi } \ \ \ \ \ in  \ \ \ R^2. \label{2.4}
\end{equation}
Let  $\Phi(S^2_i)=B_{R_i}$. It is obvious that $R_i \to +\infty$
as $i\to +\infty$.  For
$$
w_i(y)=u_i(x)+\varphi(y)=u_i(\Phi^{-1}(y))+\varphi(y),$$ we have
$$
\int_{B_{R_i}} e^{2 w_i(y)} dy= \int_{S^2_i} e^{2u_i}dx:=a_i,
$$
and \begin{align*} \int_{B_{R_i}} |\nabla w_i|^2
dy&=\int_{B_{R_i}}|\nabla (u_i\circ \Phi^{-1})|^2 dy+2
\int_{B_{R_i}}\nabla (u_i\circ \Phi^{-1}) \cdot \nabla
\varphi dy +\int_{B_{R_i}}|\nabla \varphi|^2d y\\
&= \int_{S^2_i} |\nabla u_i|^2 dx+2 \int_{S^2_i} u_i
dx+\int_{B_{R_i}}|\nabla \varphi|^2dy,
\end{align*}
where we use the fact that $ \varphi $ satisfies (\ref{2.4}).
Since $w_i(y)=\ln\frac 2{1+R_i^2}$  on $\partial B_{R_i}$, it
follows from Corollary \ref{LZcor1} that
$$
\int_{B_{R_i}} |\nabla w_i|^2 dy \ge 4 \pi (\ln \frac{a_i \cdot
(\frac {1+R_i^2}2)^2}{\pi R_i^2} +\frac {\pi R_i^2} {a_i \cdot
(\frac {1+R_i^2}2)^2} -1).
$$
Also, one can check that
$$
\int_{B_{R_i}}|\nabla \varphi|^2dy=4\pi[\ln(1+R_i^2)+\frac
1{1+R_i^2}-1].$$ We conclude
\begin{align*}
\int_{S^2_i} |\nabla u_i|^2 dx+2 \int_{S^2_i} u_i dx \ge  &4 \pi
(\ln \frac{a_i \cdot (\frac {1+R_i^2}2)^2}{\pi R_i^2} +\frac {\pi
R_i^2} {a_i \cdot (\frac {1+R_i^2}2)^2} -1)\\
&-4\pi[\ln(1+R_i^2)+\frac 1{1+R_i^2}-1]\\
=&4 \pi (\ln \frac{a_i \cdot ( {1+R_i^2})}{4\pi R_i^2} +\frac
{4\pi R_i^2} {a_i \cdot ({1+R_i^2})^2} -\frac 1{1+R_i^2}).
\end{align*}
Sending $i\to +\infty$, we have
$$
\int_{S^2} |\nabla u|^2 dx+2 \int_{S^2} u dx \ge 4 \pi ( \ln \frac
1{4\pi} \int_{S^2} e^{2u}dx).
$$


\section{Appendix}
For readers' convenience, we include geometric and analytic forms
of sharp Sobolev inequality on $S^n$ (for $n \ge 3$),  as well as
geometric and analytic forms of Onofri  inequality on $S^2$. These
are well-known to experts in the field.

\medskip

\noindent{\bf Sharp Sobolev inequality on $S^n$ (for $n\ge 3$):} \
{\it Let $(S^n, g_0)$ be the standard unit sphere in $R^{n+1}$ ($n
\ge 3$). For any $u \in H^1(S^n)$,
$$(\frac 1{\omega_n} \int_{S^n}|u|^{\frac{2n}{n-2}}dv_{g_0})^{\frac{n-2}n} \le \frac
1{\omega_n} \int_{S^n } ( u^2+\frac{4}{n(n-2)}|\nabla u|^2)
dv_{g_0}.
$$ The equality holds if and only if the scalar curvature of $u^{\frac 4{n-2}}g_0$ is
constant.}

\medskip

 If $\tilde g= \rho g$ is a conformal metric to the
background metric $g$, then the new scalar curvature $\tilde R$
under metric $\tilde g$ satisfies
\begin{equation}
\tilde R=\rho^{-1}R-(n-1)\rho^{-2} \Delta \rho-\frac 14 (n-1)(n-6)
\rho^{-3} |\nabla \rho|^2,
 \label{4.1}
\end{equation}
where $R$ is the scalar curvature under metric $g$.  If we write
$\rho=e^{2u}$, we have
\begin{equation}
\tilde R= e^{-2u}[R- (n-1)(n-2)|\nabla u|^2-2 (n-1)\Delta
u].\label{4.2}
\end{equation}
The normalized total scalar curvature under  metric $\tilde g$ is
defined by
\begin{equation}
E(\tilde g)=\frac{\int_{S^n} \tilde R dV_{\tilde g}}{(\int_{S^n}
dV_{\tilde g})^{(n-2)/n}}. \label{4.3}
\end{equation}

\noindent{\bf Geometric form of Sharp Sobolev inequality on $S^n$
(for $n\ge 3$) :} {\it Let $(S^n, g_0)$ be the standard unit
sphere in $R^{n+1}$ ($n \ge 3$). Then
$$
\inf_{\tilde g=\rho g_0} E(\tilde g)= n(n-1) \omega_n^{2/n}, $$
and the infimum is achieved if and only if $\tilde R$ (under
metric $\tilde g=\rho g_0$) is a constant.}

\medskip

For dimension $n=2$, we have

\noindent{\bf Onofri inequality on $S^2$:} {\it Let $(S^2, g_0)$
be the standard unit sphere in $R^{3}$. For any $u \in
W^{1,2}(S^2)$,
$$\ln (\frac 1{4 \pi} \int_{S^2}e^{2u}dx) \le \frac 1{4 \pi}
\int_{S^2 } (|\nabla u|^2+2u)dx. $$ The equality holds if and only
if the curvature under metric $e^{2u}g_0$ is constant.}

\medskip

Let $(M,g)$ be a smooth Riemann surface. For any conformal new
metric $g_1=e^{2u} g$, the corresponding Liouville energy is
defined by \begin{equation*}  L_g(g_1)=  \frac 14 \int_M \ln \frac
{g_1}{g} \cdot (R_{g_1} dV_{g_1}+ R_{g} dV_{g})
\end{equation*} where $R_g$ and $R_{g_1}$ are  twice the  Gaussian
curvatures $K_{g}$ and $K_{g_1}$ with respect to metrics $g$ and
$g_1$. Due to (\ref{4.2}),  the Liouville energy of metric $g_1$
can also be represented by
$$ L_g(g_1) =  \int_M  (|\nabla_g u|^2+ R_g u) dV_{g}.
$$

\noindent{\bf Geometric form of Onofri inequality on $S^2$:} {\it
Let $(S^2, g_0)$ be the standard unit sphere in $R^{3}$, and
$[g_0]_1=\{g=\rho g_0, \ \mbox{for \ some \ } \rho>0, \ \mbox{and}
\ \ \int_{S^2} dV_g=4 \pi\}$. Then
$$
\inf_{ g\in [ g_0]_1} L_{g_0}(g)=0,
$$
and the infimum is achieved if and only if $ R=2$ (under metric $
g=\rho g_0$).}

\end{document}